\documentclass {amsart}
\author{Dilip Raghavan}
\title{There is a Van Douwen MAD family}
\usepackage{amssymb, amsmath, amsthm, mathrsfs}
\newtheorem{Theorem}{Theorem}

\newtheorem{Lemma}{Lemma}
\newtheorem{Def}{Definition}
\newtheorem{Remark}{Remark}
\newtheorem{Convention}{Convention}
\newtheorem{Cor}{Corollary}
\newtheorem{conj}{Conjecture}
\renewcommand{\qedsymbol}{$\dashv$}

\newcommand{\restrict}{\upharpoonright}
\newcommand{\forallbutfin}{{\forall}^{\infty}}
\newcommand{\existsinf}{{\exists}^{\infty}}
\renewcommand{\c}{\mathfrak{c}}
\renewcommand{\[}{\left[}
\renewcommand{\]}{\right]}
\newcommand{\lc}{\left|}
\newcommand{\rc}{\right|}
\newcommand{\<}{\prec}
\DeclareMathOperator{\non}{non}
\DeclareMathOperator{\tr}{tr}
\DeclareMathOperator{\dom}{dom}
\newcommand{\Pset}{\mathcal{P}}
\renewcommand{\P}{\mathbb{P}}
\newcommand{\M}{\mathcal{M}}
\newcommand{\B}{\mathscr{B}}
\newcommand{\A}{{\mathscr{A}}}
\newcommand{\C}{{\mathscr{C}}}
\newcommand{\D}{{\mathscr{D}}}
\newcommand{\I}{{\mathcal{I}}}
\newcommand{\E}{{\mathcal{E}}}
\newcommand{\X}{{\mathscr{X}}}
\newcommand{\U}{{\mathcal{U}}}
\newcommand{\V}{\mathbf{V}}
\newcommand{\IA}{{\I}_0\left(\A\right)}
\begin{document}
\maketitle
\section{A Van Douwen MAD family in ZFC}
	\begin{Def}
		An a.d.\ family $\mathscr{A} \subset {\omega}^{\omega}$ is called a Van Douwen MAD family if for any infinite partial function $f$ from $\omega$ to $\omega$ there is $h \in \mathscr{A}$ such that $\left|f \cap h \right| = \omega$. 
	\end{Def}

It is easily seen that Van Douwen MAD families exist under CH, and more generally under MA. The question of whether they always exist was raised by E. van Douwen and A. Miller. It occurs as problem 4.2 in A. Miller's problem list~\cite{Mi2}. Zhang~\cite{Z2} discusses this problem and proves that Van Douwen MAD families of various sizes exist in certain forcing extensions.
 
In this section we will prove in ZFC that there is a Van Douwen MAD family of size Continuum. The starting point for our construction is the following well known characterization of the cardinal $\non{(\M)}$, due to Bartoszynski. The reader may consult \cite{Bar} or \cite{BJ} for a proof of this.

	\begin{Def}
		$\non{(\M)}$ is the least size of a non meager set of reals.		
	\end{Def}    

	\begin{Def}
		Let $h \in {\omega}^{\omega}$ be such that $\forall n \in \omega \left[ h(n) \geq 1 \right]$. An $h$-slalom is a function $S: \omega \rightarrow {\left[ \omega \right]}^{< \omega}$ such that for all $n \in \omega$, $\left|S(n)\right| \leq h(n)$.
	\end{Def}

	\begin{Theorem}[Bartoszynski~\cite{Bar}]\label{thm:Bar}
		Let $\kappa$ be an infinite cardinal. The following are equivalent:
			\begin{enumerate}
				\item
					Every set of reals of size less than $\kappa$ is meager.
				\item
					For every family $F \subset {\omega}^\omega$ with $\left| F \right| < \kappa$, there is an infinite partial function $g$ from $\omega$ to $\omega$ such that $\forall f \in F \left[\; \left| f \cap g \right| < \omega \right]$. 			
			
				\item
					For every $h$ and for every family of $h$-slaloms $F$ with $\left| F \right| < \kappa$, there is a $g \in {\omega}^{\omega}$ such that $\forall S \in F \; \forallbutfin n \in \omega \left[\; g(n) \notin S(n) \right]$. 
			\end{enumerate}			
	\end{Theorem} 

\begin{flushright}
	{\qedsymbol}
\end{flushright}

Our first task is to strengthen condition 3 above.

\begin{Lemma} \label{lem:onetoone}
	Let $\kappa = \non{(\M)}$ and let $F$ be a family of $h$-slaloms with $\left| F \right| < \kappa$. There is a one-to-one function $g \in {\omega}^{\omega}$ such that $\forall S \in F \; \forallbutfin n \in \omega \left[  g(n) \notin S(n) \right]$.   
\end{Lemma}

\begin{proof}
	Our proof is similar to the argument in Bartoszynski~\cite{Bar}. Write F = $\langle S_{\xi}: \xi < \lambda \rangle$, where $\lambda = \left| F \right|$. Define a new function $h'$ and a family of $h'$-slaloms as follows:
		\begin{align}	
			h'(n)   					  & = \sum_{i \leq n} h(i) \notag \\
			\forall \xi < \lambda \thickspace {S}'_{\xi}(n)   & = \bigcup_{i \leq n} S_{\xi}(i). \notag
		\end{align}

Clearly, $\langle {S}'_{\xi}: \xi < \lambda \rangle$ is a family of $h'$-slaloms. Now, for each $i \in \omega$, let $T_i: \omega \rightarrow {\left[ \omega \right]}^{< \omega}$ be defined by $T_i(n) = \{i\}$. It is clear that $\langle {S}'_{\xi}: \xi < \lambda \rangle \cup \langle T_i: i \in \omega \rangle $ is a family of fewer than $\kappa$ $h'$-slaloms. Thus by 3 of Theorem \ref{thm:Bar}, we can choose $g \in {\omega}^{\omega}$ such that the following hold: 
		\begin{enumerate}
			\item
				$\displaystyle{\forall \xi < \lambda \forallbutfin n \in \omega \left[  g(n) \notin {S}'_{\xi}(n) \right]}$
			\item
				$\displaystyle{\forall i \in \omega \forallbutfin n \in \omega \left[  g(n) \notin T_{i}(n) \right]}$.		
		\end{enumerate}
		  
Property 2 implies that $g$ takes any given value only finitely often. Thus we may choose a one-to-one infinite partial function $g' \subset g$. Let $X = \dom{(g')}$. By property 1 we obviously have that for any $\xi < \lambda$, $\forallbutfin n \in \omega \left[n \in X \implies g'(n) \notin {S}'_{\xi}(n) \right]$. Let $\langle x_n : n \in \omega \rangle$ be the increasing enumeration of $X$. For $n \in \omega$, set $g''(n) = g'(x_n)$. Since $g'$ is one-to-one, $g''$ is also one-to-one. We claim that $g''$ is the function we are looking for. Indeed, fix $\xi < \lambda$. We know that $ \exists m \in \omega \forall n \geq m \left[n \in X \implies g'(n) \notin {S}'_{\xi}(n)\right]$. We will show that $\forall n \geq m \left[g''(n) \notin S_{\xi}(n)\right]$. Suppose, for a contradiction, that $g''(n) = g'(x_n) \in S_{\xi}(n)$, for some $n \geq m$. Note that we have $m \leq n \leq x_n$. Thus, by the definition of ${S}'_{\xi}$, $S_{\xi}(n) \subset {S}'_{\xi}(x_n)$. Therefore, we have that $g'(x_n) \in {S}'_{\xi}(x_n)$. But this is a contradiction because $x_n \geq m$ and $x_n \in X$.
\end{proof}

\begin{Convention}
	In what follows we will only be concerened with $h$-slaloms for $h \equiv 2^n$. We will simply refer to these as slaloms, supressing mention of $h$.
\end{Convention}

\begin{Lemma} \label{lem:a.d.slaloms}
	Let $F = \langle S_{\xi}: \xi < \lambda \rangle$ be a family of slaloms with $\lambda < \non{(\M)}$. There is a slalom $S$ such that $\forall n \in \omega \left[ \left| S(n) \right| = 2^n \right]$ and $\forall \xi < \lambda \forallbutfin n \in \omega \left[  S(n) \cap S_{\xi}(n) = 0 \right]$.   			
\end{Lemma}               

\begin{proof}
	For all $n \in \omega$ set $l_n = 2^n - 1$ and $I_n = [l_n, l_{n+1})$. For each $\xi < \lambda$ define ${S}'_{\xi}$ by stipulating that $\forall k, n \in \omega \[{S}'_{\xi}(k) = S_{\xi}(n) \ iff \  k \in I_n \]$. We have that for all $k \in \omega$, $\left|{S}'_{\xi}(k)\right| \leq \left| S_{\xi}(n)\right| \leq 2^n$, where $k \in I_n$. But if $k \in I_n$, then $2^n \leq 2^k$ and so $\left|{S}'_{\xi}(k)\right| \leq 2^k$. Therefore, $\langle {S}'_{\xi}: \xi < \lambda \rangle$ is a family of fewer than $\non{(\M)}$ many slaloms. By applying Lemma \ref{lem:onetoone} we can find a one-to-one function $g \in {\omega}^{\omega}$ such that for every $\xi < \lambda$, $\forallbutfin k \in \omega \left[g(k) \notin {S}'_{\xi}(k)\right]$. Now define $S$ by setting $S(n) = \{g(k): k \in I_n \}$. Since $g$ is one-to-one, $\left| S(n) \right| = \left| I_n \right| = 2^n$. Fix $\xi < \lambda$. We know that $\exists m \in \omega \forall k \geq m \left[ g(k) \notin {S}'_{\xi}(k)\right]$. We claim that for any $n \geq m$, $S(n) \cap S_{\xi} (n) = 0$. Suppose to the contrary that for some $n \geq m$, $g(k) \in S_{\xi} (n)$ for some $k \in I_n$. Then since $k \in I_n$, ${S}'_{\xi}(k) = S_{\xi}(n)$, and so we get that $g(k) \in {S}'_{\xi}(k)$. But this is a contradiction because $m \leq n \leq l_n \leq k$.
\end{proof}

\begin{Lemma} \label{lem:a.d.familyinslalom}
	Let $S$ be a slalom such that $\forall n \in \omega \left[ \left| S(n) \right| = 2^n \right]$. There exists an a.d.\ family $\mathscr{A} \subset {\omega}^{\omega}$ such that $\left|\mathscr{A}\right| = \mathfrak{c}$ and for every $f \in \mathscr{A}$, $\forall n \in \omega \left[ f(n) \in S(n) \right]$.
\end{Lemma}

\begin{proof}
	Since $\left|S(n)\right| = \left|{}^n2\right|$, we can assign to each $\sigma \in {}^n2$ a unique number $k_{\sigma} \in S(n)$. Now, for each $\mu \in 2^{\omega}$, define $f_{\mu} \in {\omega}^{\omega}$ by setting $f_{\mu}(n) = k_{\mu \restrict n} \in S(n)$. Suppose $\mu \neq \nu \in 2^{\omega}$. Then there is $m \in  \omega$ such that $\mu(m) \neq \nu(m)$. So for all $n > m$, $\mu \restrict n \neq \nu \restrict n$, and so $f_{\mu}(n) = k_{\mu \restrict n} \neq k_{\nu \restrict n} = f_{\nu}(n)$. Thus $\mathscr{A} = \{f_{\mu}: \mu \in 2^{\omega} \}$ is as required. 
\end{proof}

\begin{Def}
	Let $A, B \subset {\omega}^{\omega}$ be two families of functions. We will write $A \perp B$ to mean that $\forall f \in A \forall g \in B \left[ \; \left|f \cap g\right| < \omega \right]$ 
\end{Def}

The next lemma will play an important role in our construction. The proof of this lemma will use Lemma \ref{lem:a.d.familyinslalom} and is the reason why we set out to strengthen clause (3) of Thoerem \ref{thm:Bar}.

\begin{Lemma} \label{lem:largea.d.fam}
	Let $\kappa = \non{(\M)}$. Let $F = \langle f_{\alpha}: \alpha < \kappa \rangle \subset {\omega}^{\omega}$. There is a sequence $\langle {\mathscr{A}}_{\alpha}: \alpha < \kappa \rangle$ such that following hold:
	\begin{enumerate}
		\item
			${\mathscr{A}}_{\alpha} \subset {\omega}^{\omega}$ is an a.d.\ family.
		\item
			$\left|{\mathscr{A}}_{\alpha}\right| = \mathfrak{c}$.
		\item
			for all $\beta < \alpha < \kappa$, ${\mathscr{A}}_{\alpha} \perp {\mathscr{A}}_{\beta}$
		\item
			${\mathscr{A}}_{\alpha} \perp \{f_{\beta}: \beta \leq \alpha \}$.
	\end{enumerate} 
\end{Lemma}

\begin{proof}
	We will construct the family $\langle {\mathscr{A}}_{\alpha}: \alpha < \kappa \rangle$ by induction. We will simultaneously build a family of slaloms $\langle S_{\alpha}: \alpha < \kappa \rangle$ and ensure that for all $\alpha < \kappa$, $\forall f \in {\mathscr{A}}_{\alpha} \forall n \in \omega \left[f(n) \in S_{\alpha}(n) \right]$. Fix $\alpha < \kappa$ and suppose that $\langle{\mathscr{A}}_{\beta}: \beta < \alpha\rangle$ and $\langle S_{\beta}: \beta < \alpha \rangle$ are already given to us. For each $\beta \leq \alpha$, define a slalom $T_{\beta}$ by $T_{\beta}(n) = \{f_{\beta}(n)\}$. Thus, $\{ S_{\beta}: \beta < \alpha \} \cup \{ T_{\beta}: \beta \leq \alpha \}$ is a family of fewer than $\kappa$ slaloms. So we can apply Lemma \ref{lem:a.d.slaloms} to find a slalom $S_{\alpha}$ such that the following hold:
	\begin{enumerate}
		\item[(a)]
			$\forall n \in  \omega \left[ \left| S_{\alpha}(n) \right| = 2^n \right]$ 
		\item[(b)]
			$\forall \beta < \alpha \forallbutfin n \in  \omega \left[S_{\alpha}(n) \cap S_{\beta}(n) = 0 \right]$
		\item[(c)]
			$\forall \beta \leq \alpha \forallbutfin n \in  \omega \left[S_{\alpha}(n) \cap T_{\beta}(n) =0 \right]$.
	\end{enumerate}
	Property (a) allows us to apply Lemma \ref{lem:a.d.familyinslalom} to $S_{\alpha}$ to find an a.d.\ family ${\mathscr{A}}_{\alpha} \subset {\omega}^{\omega}$ with $\left|{\mathscr{A}}_{\alpha}\right| =\mathfrak{c}$ and with the property that $\forall f \in {\mathscr{A}}_{\alpha} \forall n \in \omega \left[f(n) \in S_{\alpha}(n) \right]$. Thus ${\mathscr{A}}_{\alpha}$ satisfies requirements (1) and (2). We will check requirements (3) and (4). 
Fix $f \in {\mathscr{A}}_{\alpha}$ and $g \in {\mathscr{A}}_{\beta}$ for some $\beta < \alpha$. We know that there is $m \in \omega$ such that $\forall n \geq m \left[ S_{\alpha}(n) \cap S_{\beta}(n) =0 \right]$. Since $\forall n \in \omega \[ f(n) \in S_{\alpha} (n) \wedge g(n) \in S_{\beta}(n) \]$, it follows that $\forall n \geq m \[ f(n) \neq g(n)\]$. To verify (4), fix $f \in {\mathscr{A}}_{\alpha}$ and some $\beta \leq \alpha$. Again we know that there is $m \in \omega$ such that $\forall n \geq m \[S_{\alpha}(n) \cap \{f_{\beta} (n)\} = 0 \] $ and that $\forall n \in \omega \[f(n) \in S_{\alpha}(n) \]$. Therefore, it follows that $\forall n \geq m \[f(n) \neq f_{\beta}(n) \]$.    
\end{proof}

We are now ready to construct our Van Douwen MAD family. In order to ensure that our family is Van Douwen MAD we will introduce the notion of the trace of an a.d.\ family. The idea is that if an a.d.\ family has a "sufficiently large" trace, then it must be Van Douwen MAD.
\begin{Convention}
	By Theorem \ref{thm:Bar} there is a family $F = \langle f_{\alpha}: \alpha < \non{(\M)} \rangle \subset {\omega}^{\omega}$ such that for every infinite partial function $g$ there is an $\alpha < \non{(\M)}$ such that $\left|g \cap f_{\alpha} \right| = \omega$. For the remainder of this section let us fix such a family $F$.
\end{Convention}

\begin{Convention}
	For a countable set $X$, a MAD family of subsets of X is usually required to be an infinite family. However, in what follows we adpot the convention that a MAD family on X is simply an a.d.\ family of infinite subsets of $X$ such that every infinite subset of $X$ meets some member of the family in an infinite set. Also, since we identify functions with their graphs, it makes sense to talk of a MAD family on a function $f \in  {\omega}^{\omega}$.  	
\end{Convention}

\begin{Def}
	Let $\mathscr{A} \subset {\omega}^{\omega}$ be an a.d.\ family and let $f \in {\omega}^{\omega}$. We define $\mathscr{A} \cap f = \{ f \cap h: h \in \mathscr{A} \wedge \left| f \cap h \right| = \omega \}$. Note that this is an a.d.\ family on $f$.
\end{Def}

\begin{Def}
	Let $\mathscr{A} \subset {\omega}^{\omega}$ be an a.d.\ family. The trace of $\mathscr{A}$, written $	\tr{(\mathscr{A})}$, is $\{f \in  {\omega}^{\omega}: \mathscr{A} \cap f$ is a MAD family on f $\}$.
\end{Def}

\begin{Lemma} \label{lem:largetrace}
	Let $\mathscr{A} \subset {\omega}^{\omega}$ be an a.d.\ family such that $F \subset \tr{(\mathscr{A})}$. Then $\mathscr{A}$ is Van Douwen MAD. 
\end{Lemma}

\begin{proof}
	Indeed, let $g$ be an infinite partial function. By the definition of $F$, there is $\alpha < \non{(\M)}$ such that $\left| g \cap f_{\alpha} \right| = \omega$. Since $F \subset \tr{(\mathscr{A})}$, $\mathscr{A} \cap f_{\alpha}$ is a MAD family on $f_{\alpha}$. So there is $h \in \mathscr{A}$ such that $h \cap f_{\alpha}$ meets $g \cap f_{\alpha}$ in an infinite set, whence we get that $\left| h \cap g \right| = \omega$.  
\end{proof}

\begin{Theorem} \label{thm:VMADsexist}
	There is a Van Douwen MAD family of size $\mathfrak{c}$.
\end{Theorem}

\begin{proof}
	In view of Lemma \ref{lem:largetrace}, it is enough to construct an a.d.\ family $\mathscr{A}$ of size $\mathfrak{c}$ such that $F \subset \tr{(\mathscr{A})}$. We will use Lemma \ref{lem:largea.d.fam} to do this. Fix a sequence $\langle {\mathscr{A}}_{\alpha}: \alpha < \non{(\M)} \rangle$ as in Lemma \ref{lem:largea.d.fam}. $\mathscr{A}$ will be constructed as the union of an increasing sequence of a.d.\ families. Thus, we will construct a sequence $\langle {\mathscr{C}}_{\alpha}:  \alpha < \non{(\M)} \rangle$ such that:
	\begin{enumerate}
		\item
			${\mathscr{C}}_{\alpha} \subset {\omega}^{\omega}$ is an a.d.\ family
		\item
			$\forall \beta < \alpha < \non{(\M)} \[{\mathscr{C}}_{\beta} \subset {\mathscr{C}}_{\alpha} \]$
		\item
			$f_{\alpha} \in \tr{({\mathscr{C}}_{\alpha})}$
		\item
			$\forall h \in {\mathscr{C}}_{\alpha} \exists \beta \leq \alpha \exists g \in {\mathscr{A}}_{\beta} \exists X \in  {\[\omega\]}^{\omega} \[h = f_{\beta} \restrict X \cup g \restrict {\omega \setminus X} \]$
		\item
			$\left| {\mathscr{C}}_{0}\right| = \mathfrak{c}$.
	\end{enumerate}     

	To construct ${\mathscr{C}}_{0}$, we fix a MAD family $\{a_{\xi}: \xi < \mathfrak{c}\}$ on $\omega$. Put ${\mathscr{A}}_{0} = \{g_{\xi}: \xi < \mathfrak{c}\}$. For each $\xi < \mathfrak{c}$, let $h_{\xi} = f_{0} \restrict a_{\xi} \cup g_{\xi} \restrict {{\omega} \setminus a_{\xi}} $, and put ${\mathscr{C}}_{0} = \{h_{\xi}: \xi < \mathfrak{c} \}$. We will check that ${\mathscr{C}}_{0}$ is a.d. Indeed, if $\xi_{0} < \xi_{1}$, then since $a_{\xi_{0}} \cap  a_{\xi_{1}}$ is finite, $\left| f_{0} \restrict a_{\xi_{0}} \cap f_{0} \restrict a_{\xi_{1}} \right| < \omega$. Next, since ${\mathscr{A}}_{0} \perp \{ f_{0} \}$, we have that both $f_{0} \restrict a_{\xi_{0}} \cap g_{\xi_{1}} \restrict {{\omega} \setminus a_{\xi_{1}} }$ and $f_{0} \restrict a_{\xi_{1}} \cap g_{\xi_{0}} \restrict {{\omega} \setminus a_{\xi_{0}} }$ are finite. Finally, since ${\mathscr{A}}_{0}$ is an a.d.\ family, we know that $\left| g_{\xi_{0}} \restrict {{\omega} \setminus a_{\xi_{0}} } \cap g_{\xi_{1}} \restrict {{\omega} \setminus a_{\xi_{1}} } \right| < \omega$. Thus, we conclude that $\left| h_{\xi_{0}} \cap h_{\xi_{1}} \right| < \omega$. Next, it is clear from the construction that $f_{0} \in \tr{({\mathscr{C}}_{0})}$, and that ${\mathscr{C}}_{0}$ satisfies clauses (4) and (5). 

	To continue the construction, suppose that we are given the sequence $\langle {\mathscr{C}}_{\beta}:  \beta < \alpha \rangle$. Set $\mathscr{C} = \bigcup {\mathscr{C}}_{\beta}$ and consider $\mathscr{C} \cap f_{\alpha}$. This is an a.d.\ family on $f_{\alpha}$. If it is a MAD family (either finite or infinite), then $f_{\alpha}$ is already in $\tr{(\mathscr{C})}$, and there is nothing more to be done. In this case, we set ${\mathscr{C}}_{\alpha} = \mathscr{C}$. So, say that $\mathscr{C} \cap f_{\alpha}$ is not MAD. We can extend it to a MAD family, say $\mathscr{B}$, on $f_{\alpha}$. Consider the family $\{Y \in {\[\omega\]}^{\omega}: f_{\alpha} \restrict  Y \in \mathscr{B} \setminus \left( \mathscr{C} \cap f_{\alpha}\right) \}$. Note that this is an a.d.\ family on $\omega$. We may assume WLOG that it has size $\c$. Let $\{a_{\xi}: \xi < \mathfrak{c} \}$ enumerate this family. Put ${\mathscr{A}}_{\alpha} = \{g_{\xi}: \xi < \mathfrak{c}\}$.  For each $\xi < \mathfrak{c}$ set $h_{\xi} = f_{\alpha} \restrict a_{\xi} \cup g_{\xi} \restrict {\omega \setminus a_{\xi}}$, and put $\mathscr{D} = \{h_{\xi}: \xi < \mathfrak{c} \}$. It is easily argued, as for ${\mathscr{C}}_{0}$, that $\mathscr{D}$ is a.d. We will check that $\mathscr{C} \perp \mathscr{D}$. Fix $h \in \mathscr{C}$ and $\xi < \mathfrak{c}$. If $h \cap f_{\alpha}$ is finite, then so is $h \cap f_{\alpha} \restrict a_{\xi}$. On the other hand, if $h \cap f_{\alpha}$ is infinite, then $h \cap f_{\alpha} \in \mathscr{C} \cap f_{\alpha}$. But then $\left| f_{\alpha} \restrict a_{\xi} \cap h \right| < \omega$ because $\mathscr{B}$ is an a.d.\ family. Thus in either case, $h \cap f_{\alpha} \restrict a_{\xi}$ is finite. To deal with $h \cap g_{\xi} \restrict {\omega \setminus a_{\xi}}$, by clause (4), we know that for some $\gamma \leq \beta < \alpha$, $h = f_{\gamma} \restrict X \cup g \restrict {\omega \setminus X}$, where $X \in {\[\omega\]}^{\omega}$ and $g \in {\mathscr{A}}_{\gamma}$. But since ${\mathscr{A}}_{\alpha} \perp \{f_{\gamma} \}$, $\left|f_{\gamma} \restrict X \cap g_{\xi} \restrict {\omega \setminus a_{\xi}} \right| < \omega$, and since ${\mathscr{A}}_{\alpha} \perp {\mathscr{A}}_{\gamma}$, we know that $\left| g_{\xi} \restrict {\omega \setminus a_{\xi}} \cap  g \restrict {\omega \setminus X} \right|< \omega$. Therefore, $h \cap g_{\xi} \restrict {\omega \setminus a_{\xi}}$ is also finite, and so $\left|h \cap h_{\xi} \right| < \omega$. Hence, we can define ${\mathscr{C}}_{\alpha} = \mathscr{C} \cup \mathscr{D}$.

	Now, it is clear that ${\mathscr{C}}_{\alpha}$ satisfies clauses (1), (2) and (4). We just need to verify that $f_{\alpha} \in \tr{({\mathscr{C}}_{\alpha})}$. So we need to check that ${\mathscr{C}}_{\alpha} \cap f_{\alpha}$ is a MAD family on $f_{\alpha}$. But clearly ${\mathscr{C}}_{\alpha} \cap f_{\alpha} = \mathscr{C} \cap f_{\alpha} \cup \mathscr{D} \cap f_{\alpha}$. Fix $X \in {\[\omega\]}^{\omega}$. Since $\mathscr{B}$ is a MAD family on $f_{\alpha}$, there is $Y \in {\[\omega\]}^{\omega}$ such that $f_{\alpha} \restrict Y \in \mathscr{B}$ and $\left| f_{\alpha} \restrict X \cap f_{\alpha} \restrict Y \right| = \omega$. If $f_{\alpha} \restrict Y \in \mathscr{C} \cap f_{\alpha}$, then we are done. If it is not, then $Y = a_{\xi}$ for some $\xi < \mathfrak{c}$. It follows that $\left|f_{\alpha} \restrict X  \cap h_{\xi} \right|= \omega$. But since $ h_{\xi} \in \mathscr{D}$, we are done.     
\end{proof}
\section{Definability of MAD families in ${\omega}^{\omega}$}
	Our next task is to investigate the definability of a.d.\ families in ${\omega}^{\omega}$. We will first prove that if $\A$ is an analytic MAD family in ${\omega}^{\omega}$, then $\A$ must satisfy some strong constraints (Theorem \ref{thm:trivialtrace}). This will immediately imply that Van Douwen MAD families cannot be analytic. We will then show that this is a strengthening of a result of J. Steprans~\cite{KSZ} that strongly MAD families cannot be analytic. Next, we will show that it is consistent to have MAD families in ${\omega}^{\omega}$ that satisfy these strong constraints. Finally, we will argue that 
analytic MAD families cannot satisfy these constraints if they have some additional combinatorial properties.
\begin{Def}
	Let $\A \subset {\[{\omega} \times {\omega} \]}^{\omega}$ be an a.d.\ family and let $X \in \[{\omega} \times {\omega}\]^{\omega}$. We say that $X$ avoids $\A$ if for any finite collection $\{h_0, \dotsc, h_n\} \subset \A$, $\left| X \setminus {h_0 \cup \dotsb \cup h_n} \right| = \omega$.
\end{Def}

\begin{Theorem} \label{thm:trivialtrace}
	Let $\A \subset {\omega}^{\omega}$ be an a.d.\ family and let $X \in \[{\omega} \times {\omega}\]^{\omega}$ avoid $\A$. Suppose that $\A$ is analytic in $ {\omega}^{\omega}$. There is $Y \in \[ X \]^{\omega}$ such that $\forall h \in \A \[\; \left| h \cap Y \right|  < \omega \]$.
\end{Theorem}

\begin{proof}
	Let us give the space $2^{X}$ the Tychonoff product topology, with $2$ having the discrete topology. Since $X$ is a countable set, this is homeomorphic to $2^{\omega}$ with the usual topology. Define a map $\Psi: {\omega}^{\omega} \rightarrow 2^{X}$ by stipulating that $\forall \langle n, m \rangle \in X \[ \Psi(f)\left( \langle n, m \rangle \right) = 1 \leftrightarrow \langle n, m \rangle \in f \]$. Thus $\Psi(f)$ is the characteristic function of $X \cap f$. We will argue that this map is continious. Fix finitely many members $\langle n_0, m_0 \rangle, \dotsc, \langle n_k, m_k \rangle \in X$ and $\langle n^{0}, m^{0} \rangle, \dotsc, \langle n^l, m^l \rangle \in X$. A basic open subset of $2^X$ is of the form $ U = \{\chi \in 2^X: \chi \left(\langle n_i, m_i \rangle\right) = 0 \; \forall i \leq k \wedge \chi \left(\langle n^i, m^i \rangle \right) = 1 \; \forall i \leq l  \}$. Thus ${\Psi}^{-1} (U) = \{f \in {\omega}^{\omega}: f(n_i) \neq m_i \; \forall i \leq k \wedge f(n^i) = m^i \; \forall i \leq l \}$. It is clear that this is an open subset of ${\omega}^{\omega}$. Hence $\Psi$ is continious. Therefore, ${\Psi}''{\A}$ is an analytic subset of $2^X$. It is the set of characteristic functions of elements of $\{ h \cap X: h \in \A \}$. We are only interested in the infinite elements of this set. So we will put $ \B =  {\Psi}''{\A} \cap \{ \chi \in 2^X: \existsinf \langle n, m \rangle \in X \[\chi \left( \langle n, m \rangle \right) = 1 \] \}$. It is clear that $\B$ is also analytic. $\B$ is the set of characteristic functions of elements of $\A \cap X = \{h \cap X: h \in \A \wedge \left|h \cap X \right| = \omega \}$. Now, $\A \cap X$ is an a.d.\ family on $X$. By a theorem of Mathias~\cite{M} we know that there are no analytic MAD families on $X$. Therefore, if $\A \cap X$ is infinite, it is not MAD on $X$, and we will get the conclusion of the theorem. On the other hand, if $\A \cap X$ is finite, then since $X$ avoids $\A$, $Y = X \setminus \bigcup \left(\A \cap X \right)$ will satisfy the conclusion of the theorem. Hence, either way, the theorem is proved.     
\end{proof}

\begin{Def}
	An a.d.\ family $\A \subset {\omega}^{\omega}$ is said to have trivial trace if no member of $\tr{(\A)}$ avoids $\A$.
\end{Def}

\begin{Cor} \label{cor:trivialtrace}
	Suppose $\A \subset {\omega}^{\omega}$ is an analytic a.d.\ family. Then $\A$ has trivial trace.	
\end{Cor}

\begin{proof}
	If $f$ is a member of $\tr{(\A)}$ which avoids $\A$, then putting $f = X$ in Theorem \ref{thm:trivialtrace} will give a contradiction.
\end{proof}

\begin{Cor} \label{cor:noanalyticVMADs}
	There are no analytic Van Douwen MAD families in ${\omega}^{\omega}$.
\begin{flushright}
	\qedsymbol
\end{flushright}
\end{Cor}

	Juris Steprans~\cite{KSZ} introduced the notion of a strongly MAD family and proved that they can't be analytic.
\begin{Def}
	A MAD family $\A \subset {\omega}^{\omega}$ is strongly MAD if for every collection $\{f_i: i \in \omega \} \subset {\omega}^{\omega}$ where each $f_i$ avoids $\A$, there is $h \in \A$ such that $\forall i \in \omega \[ \; \left| f_i \cap h\right| = \omega \]$.
\end{Def}

\begin{Lemma} \label{lem:SMADimpliesstronglyVMAD}
	Let $\A \subset {\omega}^{\omega}$ be strongly MAD. Let $\{g_i: i \in \omega\}$ be a collection of infinite partial functions from $\omega$ to $\omega$ such that each $g_i$ avoids $\A$. There is $h \in \A$ such that $\forall i \in \omega \[\; \left| h \cap g_i\right| = \omega \]$. In particular, strongly MAD families are Van Douwen MAD.
\end{Lemma}

\begin{proof}
	Let $h_0 \neq h_1$ be two distinct members of $\A$. For each $i \in \omega$, let $a_i = \dom {(g_i)}$ and let $b_i = {\omega} \setminus {a_i}$. For each $i \in \omega$, define $f_{i}^{0} = g_i \cup h_0 \restrict b_i$ and $f_{i}^{1} = g_i \cup h_1 \restrict b_i$. Since $g_i$ avoids $\A$, both $f_{i}^{0}$ and $f_{i}^{1}$ avoid $\A$. Thus $\{f_{i}^{j}: i \in \omega \wedge j \in 2 \}$ is a countable collection of total functions avoiding $\A$. So we may choose $h \in \A$ such that $\forall i \in \omega \forall j \in 2 \[\; \left| h \cap f_{i}^{j} \right| = \omega \]$. We will show that $\forall i \in \omega \[\; \left| g_i \cap h \right| = \omega \]$. If $g_i \cap h$ is finite, then since both $f_{i}^{0} \cap h$ and $f_{i}^{1} \cap h$ are infinite, it follows that $\left| h_0 \cap h \right| = \omega$ and that $\left|h_1 \cap h \right| = \omega$. But since $\A$ is an a.d.\ family this means that $h = h_0$ and $h = h_1$, which is a contradiction. 
\end{proof}
\begin{Cor}[Steprans~\cite{KSZ}]
	There are no analytic strongly MAD families in ${\omega}^{\omega}$. \begin{flushright} \qedsymbol \end{flushright}
\end{Cor}
\begin{Remark}
	Corollary 2 is strictly stronger than Corollary 3. It is easy to modify the construction in Theorem 2 to ensure that the Van Douwen MAD family constructed there is not strongly MAD.
\end{Remark}

	It is an open problem whether there are any analytic MAD families in ${\omega}^{\omega}$. In fact, it is not even known if a MAD family in ${\omega}^{\omega}$ can be closed. Since Theorem \ref{thm:trivialtrace} puts a strong restriction on such MAD families, one might conjecture that there are no MAD families that satisfy the conclusion of Theorem \ref{thm:trivialtrace} at all. However, we will show below that this is consistently false. We will first argue that it is sufficient to build a MAD family with trivial trace.  
\begin{Lemma} \label{lem:trivialtraceistrivialtrace}
	Let $\A \subset {\omega}^{\omega}$ be a MAD family with trivial trace. Suppose $X \in \[{\omega} \times {\omega} \]^{\omega}$ avoids $\A$. There is $Y \in \[X \]^{\omega}$ such that $\forall h \in \A \[ \; \left|h \cap Y \right| < \omega \]$. 
\end{Lemma}

\begin{proof}
	Let $\A \cap X = \{ h \cap X: h \in \A \wedge \lc h \cap X \rc = \omega \}$. If $\A \cap X$ is finite, then since $X$ avoids $\A$, $Y = X \setminus \bigcup \left( \A \cap X\right)$ will be as required. So assume that $\A \cap X$ is infinite. Choose a countably infinite collection $\{h_i: i \in \omega \} \subset \A$ such that $\lc h_i \cap X \rc = \omega$ for each $i$, and put $p_i = h_i \cap X$. Thus $\{p_i: i \in \omega \}$ forms an a.d family of infinite partial functions. We may choose infinite partial functions $g_i \subset p_i$ such that $\forall i < j < \omega \[\dom{(g_i)} \cap \dom{(g_j)} = 0  \]$. Now if we put $g = \bigcup g_i$, then $g$ is an infinite partial function and $g \subset X$. Since $g$ has infinite intersection with infinitely many things in $\A$, it is clear that $g$ avoids $\A$. Let $a = \dom{(g)}$ and let $b = {\omega} \setminus a$. Choose $h \in \A$ and put $f = g \cup h \restrict b$. Obviously, $f$ is a total function avoiding $\A$. So $f \notin \tr {(\A)}$. Therefore, we may choose an infinite partial function $p \subset f$ such that $\forall h \in \A \[\lc h \cap p \rc < \omega \]$. Clearly, since $ \lc p \cap h \restrict b \rc < \omega$, we have that $\lc p \cap g \rc = \omega$. Thus, $Y = p \cap g$ is as required. 
\end{proof}
\begin{Def}
	Let $\B$ be a family of infinite partial functions from $\omega$ to $\omega$. We will say that $\B$ has an a.d.\ base if there is an a.d.\ family $\C \subset \B$ such that $\forall f \in \B \exists g \in \C \[ f \subset g\]$.
\end{Def}
\begin{Lemma} \label{lem:aisbigimpliesnotMAD}
	Assume that $\mathfrak{a} = \mathfrak{a_{\mathfrak{e}}} = \mathfrak{c}$. Let $\B$ be a family of infinite partial functions with an a.d.\ base, and suppose that $\lc \B \rc < \c$. Let $f \in {\omega}^{\omega}$ avoid $\B$. There is $h \in {\omega}^{\omega}$ such that:
	\begin{enumerate}
		\item
			$\forall g \in \B \[\; \lc h \cap g  \rc < \omega \]$
		\item
			$\lc h \cap f \rc = \omega$.
	\end{enumerate}    
\end{Lemma}
\begin{proof}
	Let $\C = \{g_{\alpha}: \alpha < \kappa \} \subset \B$ be an a.d.\ base for $\B$, with $\kappa < \c$. Consider $\C \cap f = \{g_{\alpha} \cap f: \alpha < \kappa \wedge \lc g_{\alpha} \cap f \rc = \omega \}$. This is an a.d.\ family on $f$. Since $f$ avoids $\B$, $\mathfrak{a} = \c$ and $\kappa < \c$, $\C \cap f$ cannot be a MAD family (either finite or infinite) on $f$. So there is an infinite partial function $p \subset f$ such that $\forall \alpha < \kappa \[\lc g_{\alpha} \cap p \rc < \omega \]$. Note that this means that $\forall g \in \B \[\lc g \cap p \rc < \omega \]$. Let $a = \dom{(p)}$ and let $b = {\omega} \setminus {a}$. Now, let $\D = \{ g_{\alpha} \restrict b: \alpha < \kappa \wedge \lc g_{\alpha} \restrict b \rc =\omega \}$. This is an a.d.\ family of functions in ${\omega}^{b}$ of size atmost $\kappa$. Since $\mathfrak{a_{\mathfrak{e}}} = \mathfrak{c}$, $\D$ cannot be maximal in ${\omega}^{b}$. So there is a function $q: b \rightarrow {\omega}$ such that $\forall \alpha < \kappa \[ \lc g_{\alpha} \cap q \rc < \omega \]$. Once again, this means that $\forall g \in \B \[\lc q \cap g \rc < \omega \]$. Now, it is clear that $h = p \cup q$ is as required. 
\end{proof}
\begin{Theorem} \label{thm:aisbigimpliesexiststrivialtraceMAD}
	Assume $\mathfrak{a} = \mathfrak{c}$. There is a MAD family $\A \subset {\omega}^{\omega}$ with trivial trace.	
\end{Theorem}
\begin{proof}
	Suppose first that $\mathfrak{a_{\mathfrak{e}}} < \mathfrak{a}$ (it is unknown if this situation is consistent). Then any MAD family $\A \subset {\omega}^{\omega}$ of size $\mathfrak{a_{\mathfrak{e}}}$ will have trivial trace because for any $f \in {\omega}^{\omega}$, $\lc \A \cap f \rc < \mathfrak{a} $. So we may assume that $\mathfrak{a} = \mathfrak{a_{\mathfrak{e}}} = \mathfrak{c}$. Let $\langle f_{\alpha}: \alpha < \c \rangle$ enumerate ${\omega}^{\omega}$. A non-principal, proper ideal on $\omega$ is said to be dense if $\forall a \in {\[{\omega}\]}^{\omega} \exists b \in {\[a\]}^{\omega} \[ b \in \I \]$. Note that such ideals always exist. For the rest of the proof, let us fix one such ideal, $\I$. We will construct the MAD family $\A$ by induction, as the union of an increasing sequence of a.d.\ families. In fact, we will build three sequences $\langle {\A}_{\alpha}: \alpha < \c  \rangle$, $\langle {\B}_{\alpha}: \alpha < \c \rangle$ and $\langle {\C}_{\alpha}: \alpha < \c \rangle$ such that the following hold:
	\begin{enumerate}
		\item
			${{\A}_{\alpha}} \subset {\omega}^{\omega}$ is an a.d family, with $\lc {\A}_{\alpha} \rc \leq \lc \alpha \rc$ 
		\item
			${\B}_{\alpha}$ is a family of infinite partial functions, with $\lc {\B}_{\alpha} \rc \leq \lc \alpha \rc$
		\item
			${\C}_{\alpha} \subset {\B}_{\alpha}$ is an a.d.\ base for ${\B}_{\alpha}$ 
		\item	
			$\forall \alpha < \beta < \c \[ {\A}_{\alpha} \subset {\A}_{\beta} \wedge {\B}_{\alpha} \subset {\B}_{\beta} \wedge {\C}_{\alpha} \subset {\C}_{\beta} \]$
		\item
			$\forall g \in {\B}_{\alpha} \[\dom{(g)} \in \I \]$
		\item
			$\forall h \in {\A}_{\alpha} \forall g \in {\B}_{\alpha} \[\; \lc h \cap g \rc < \omega \]$.
		\item
			if $f_{\alpha}$ avoids $\bigcup \{{\A}_{\beta}: \beta < \alpha \}$, then there is $g \in {\B}_{\alpha}$ so that $g \subset f_{\alpha}$	
		\item
			if $f_{\alpha}$ is a.d.\ from $\bigcup \{{\A}_{\beta}: \beta < \alpha \}$, there is $h \in {\A}_{\alpha}$ so that $\lc h \cap f_{\alpha} \rc = \omega$.
	\end{enumerate}

	$\A$ will be $\bigcup {\A}_{\alpha}$. Clauses (1) and (8) ensure that $\A$ is a MAD family in ${\omega}^{\omega}$. Clauses (6) and (7) ensure that $\A$ has trivial trace. It is easy to see that clause (5) is necessary becasue if $\A$ is a MAD family with trivial trace, then $\{a \in {\[{\omega}\]}^{\omega}: \exists p \in {\omega}^{a} \[p \ \text{is a.d.\ from} \ \A \] \}$ is a proper dense ideal on $\omega$. Finally, clauses (2) and (3) will allow us to continue the construction just from the assumption $\mathfrak{a} = \mathfrak{a_{\mathfrak{e}}} = \mathfrak{c}$.

	Fix $\alpha < \c$ and suppose that $\langle {\A}_{\beta}: \beta < \alpha  \rangle$, $\langle {\B}_{\beta}: \beta < \alpha \rangle$ and $\langle {\C}_{\beta}: \beta < \alpha \rangle$ are given to us. Set $\B = \bigcup {\B}_{\beta}$, $\C = \bigcup {\C}_{\beta}$ and $\D = \bigcup {\A}_{\beta}$. Note that $\C$ is an a.d.\ base for $\B$. If $f_{\alpha}$ does not avoid $\D$, then there is nothing to be done. In this case, we simply set ${\A}_{\alpha} = \D$, ${\B}_{\alpha} = \B$ and ${\C}_{\alpha} =\C$. From now on, let us assume that $f_{\alpha}$ avoids $\D$. Suppose there is a $g \in \B$ such that $\lc g \cap f_{\alpha} \rc = \omega$. Since $\I$ is a dense ideal, we can find an infinite partial function $g_{0} \subset g \cap f_{\alpha}$, with $\dom{(g_{0})} \in \I$. We set ${\B}_{\alpha} = \B \cup \{g_{0}\}$ and ${\C}_{\alpha}= \C$. It is clear that $\C$ is still an a.d.\ base for ${\B}_{\alpha}$. Moreover, $g_{0}$ is a.d.\ from $\D$. On the other hand, if $f_{\alpha}$ is a.d.\ from $\B$, we can proceed as follows. Consider $\D \cap f_{\alpha}$. This is an a.d.\ family on $f_{\alpha}$. Since $\lc \D \rc < \c$, we can find an infinite partial function $p \subset f_{\alpha}$ so that $\forall h \in \D \[\lc p \cap h \rc < \omega \]$. Since $\I$ is a dense ideal, there is an infinite partial function $g_{1} \subset p$ with $\dom{(g_{1})} \in \I$. Now, we define ${\B}_{\alpha} = \B \cup \{g_{1}\}$ and ${\C}_{\alpha} = \C \cup \{g_{1} \}$. Note that because of our assumption that $f_{\alpha}$ is a.d.\ from $\B$, $g_{1}$ is a.d.\ from $\C$. Thus ${\C}_{\alpha}$ is an a.d.\ base for ${\B}_{\alpha}$. Also, by our choice of $p$, we have that $\forall h \in \D \[ \lc h \cap g_{1} \rc < \omega \] $. This completes the definition of ${\B}_{\alpha}$ and ${\C}_{\alpha}$. We now define ${\A}_{\alpha}$. Again, we will proceed by cases. Suppose that $f_{\alpha}$ is not a.d.\ from $\D$. In this case, we may set ${\A}_{\alpha} = \D$. Note that we have already ensured above that everything in ${\B}_{\alpha}$ is a.d from $\D$. So clause (6) will be satisfied. All the other clauses are immediate. Now, let us consider the case when $f_{\alpha}$ is a.d.\ from $\D$. It is clear that ${\C}_{\alpha} \cup \D$ is an a.d.\ base for ${\B}_{\alpha} \cup \D$. Also, $\lc {\B}_{\alpha} \cup \D \rc < \c$. We claim that $f_{\alpha}$ avoids ${\B}_{\alpha} \cup \D$. For suppose that there are finitely many $g^{0}, \dotsc, g^{i} \in {\B}_{\alpha}$ and $h^{0}, \dotsc, h^{j} \in \D$ such that $f_{\alpha} {\subset}^{\ast} g^{0} \cup \dotsb \cup g^{i} \cup h^{0} \cup \dotsb \cup h^{j}$. Since $f_{\alpha}$ is a.d.\ from $\D$, it must be the case that $f_{\alpha} {\subset}^{\ast} g^{0} \cup \dotsb \cup g^{i}$. But $f_{\alpha}$ is a total function. So we must have that $\omega \: {\subset}^{\ast} \dom{(g^{0})} \cup \dotsb \cup \dom{(g^{i})}$. However, we know that $\I$ is a proper ideal and that $\dom{(g^{l})} \in \I$ for all $0 \leq l \leq i$. So this is impossible. Therefore, $f_{\alpha}$ must avoid ${\B}_{\alpha} \cup \D$. So we can apply Lemma \ref{lem:aisbigimpliesnotMAD} to find $h \in {\omega}^{\omega}$ which is a.d.\ from ${\B}_{\alpha} \cup \D$ and so that $\lc h \cap f_{\alpha} \rc = \omega$. Now, we can set ${\A}_{\alpha} = {\D} \cup \{h\}$. It is easy to see that clauses $(1)-(8)$ are all satisfied, and so we are done.
\end{proof}
	We do not know if this construction can be carried out in ZFC. But we conjecture that it cannot.
	\begin{conj}
		It is consistent with ZFC that every MAD family in ${\omega}^{\omega}$ has a non trivial trace.
	\end{conj}
	Theorem \ref{thm:aisbigimpliesexiststrivialtraceMAD} implies that it is consistent to have a MAD family with trivial trace. However, it may still be the case that analytic MAD families cannot have trivial trace. We will investigate this possibility next. We will show that analytic MAD satisfying certain extra combinatorial properties cannot have trivial trace, and hence, cannot exist. We will use a partition theorem proved by Taylor~\cite{BT} and extended by Blass~\cite{Bl}.
\begin{Theorem}[Taylor, see \cite{Bl} Theorem 4] \label{thm:blass}
	Let $\U$ be a P-point on $\omega$ and let $\X \subset {\[\omega\]}^{\omega}$ be an analytic set. There is a set $E \in \U$ and a function $f \in {\omega}^{\omega}$ such that $\X$ contains all or none of the infinite subsets $F$ of $E$ that satisfy
	\begin{equation*}
		\tag{$\ast$} \forall i, j \in F \[i < j \implies f(i) < j \].
	\end{equation*}
\end{Theorem} \begin{flushright} \qedsymbol \end{flushright}
\begin{Convention}
	We will apply Theorem \ref{thm:blass} to an $\U$ on $\omega \times \omega$ and an $\X \subset {\[\omega \times \omega\]}^{\omega}$. In order to make sense of the condition $\left(\ast\right)$ in Theorem \ref{thm:blass}, we must have a well ordering of $\omega \times \omega$ in type $\omega$. Let us arbitrarily choose such an ordering $\<$.
\end{Convention}
\begin{Lemma} \label{lem:nohomogeneoussetavoids}
	Let $\A \subset {\omega}^{\omega}$ be an analytic a.d.\ family. Let $E \in {\[\omega \times \omega\]}^{\omega}$ be a set such that $\existsinf h \in \A\[\lc h \cap E \rc = \omega \]$. Let $\X = \{F \in {\[\omega \times \omega\]}^{\omega}: \exists h \in \A \[\lc h \cap F \rc = \omega\] \}$. Let $f \in {\left(\omega \times \omega \right)}^{\omega}$. There are infinite sets $F_0$ and $F_1$ in ${\[E\]}^{\omega}$ such that $F_0 \in \X$, $F_1 \notin \X$ and
	\begin{align*}
		\tag{${\ast}_0$} & \forall \langle i_0, j_0 \rangle, \langle i_1, j_1 \rangle \in F_0 \[\langle i_0, j_0 \rangle \: \< \langle i_1, j_1 \rangle \implies f(\langle i_0, j_0 \rangle) \< \langle i_1, j_1 \rangle \] \\
	\tag{${\ast}_1$} & \forall \langle k_0, l_0 \rangle, \langle k_1, l_1 \rangle \in F_1 \[\langle k_0, l_0 \rangle \< \langle k_1, l_1 \rangle \implies f(\langle k_0, l_0 \rangle) \< \langle k_1, l_1 \rangle \].
	\end{align*}
\end{Lemma}
\begin{proof}
	Choose $h \in \A$ such that $\lc h \cap E \rc = \omega$. We may choose, by recursion, an infinite set $F_0 \subset h \cap E$ that satisfies $\left( {\ast}_0 \right)$ above. It is clear that $\lc F_0 \cap h \rc = \omega$, and hence that $F_0 \in \X$. To get $F_1$, we will use Theorem \ref{thm:trivialtrace}. Note that $E$ avoids $\A$. So there is $F \in {\[E\]}^{\omega}$ such that $F$ is a.d.\ from $\A$. Once again, we may choose, by recursion, an infinite set $F_1 \subset F$ that satisfies $\left( {\ast}_1 \right)$ above. It is clear that $F_1$ is a.d.\ from $\A$, and hence that $F_1 \notin \X$.
\end{proof}
\begin{Def}
	Let $A$ be a countable set and let $\I$ be a non-principal ideal on $A$. Let $\E = {\[A\]}^{\omega} \setminus \I$. We say that $\E$ is a P-coideal on $A$ if whenever $E_0 \supset E_1 \supset \dotsb$ is a sequence of sets in $\E$, there a set $E \in \E$ such that $\forall n \in \omega \[E {\subset}^{\ast} E_n \]$. 
\end{Def}
\begin{Theorem} \label{thm:IisnotPcoideal}
	Let $\A \subset {\omega}^{\omega}$ be an a.d.\ family. Let $\X = \{F \in {\[\omega \times \omega\]}^{\omega}: \exists h \in \A \[\lc h \cap F \rc = \omega\] \}$ and let ${\E}_0 = \{E \in {\[\omega \times \omega\]}^{\omega}: \existsinf h \in \A \[ \lc h \cap E \rc = \omega \] \}$. If there is a P-coideal $\E$ on $\omega \times \omega$ with $\E \subset {\E}_0$, then $\A$ is not analytic.  
\end{Theorem}
\begin{proof}
	By definition, there is a non-principal ideal $\I$ such that $\E = {\[\omega \times \omega\]}^{\omega} \setminus \I$. Let $\P$ be the forcing notion $\Pset \left( \omega \times \omega \right)\slash \I$. Since $\E$ is a P-coideal, $\P$ is countably closed and hence does not add any reals. Moreover, $\P$ generically adds a P-point $\U \subset \E$. Now, suppose for a contradiction that $\A$ is analytic. Identifying ${\omega}^{\omega}$ with a $G_{\delta}$ subset of $\Pset \left( \omega \times \omega \right)$ in the natural way makes $\A$ into an analytic subset of $\Pset \left( \omega \times \omega \right)$. This implies that $\X$ is analytic because it has a ${\Sigma}^{1}_{1}$ defintion. As $\P$ does not add any reals, $\X $ is still an analytic set in $\V\[\U\]$ with the same defintion. Now, in $\V\[\U\]$, we may apply Theorem \ref{thm:blass} to find a set $E \in \U$ and a function $f \in {\left( \omega \times \omega \right)}^{\omega}$ such that $\X$ contains all or none of the infinite subsets $F$ of $E$ that satisfy
	\begin{equation*}
	\tag{$\ast$} \forall \langle i, j\rangle, \langle k, l\rangle \in F \[\langle i, j\rangle \< \langle k, l\rangle \implies f\left(\langle i, j\rangle \right) \< \langle k, l \rangle \].
	\end{equation*}
	But $\P$ does not add any reals. Therefore, $E$ and $f$ are in the ground model $\V$. Note that $E \in \E \subset {\E}_0$ because $\U \subset \E$. This allows us to apply Lemma \ref{lem:nohomogeneoussetavoids} in $\V$ to find $F_0, F_1 \in {\[E\]}^{\omega}$ satisfying $\left({\ast}_0\right)$ and $\left( {\ast}_1 \right)$ of Lemma \ref{lem:nohomogeneoussetavoids} with $F_0 \in \X$ and $F_1 \notin \X$. But now, $F_0, F_1 \in \V \[\U\]$ still satisfy $\left({\ast}_0\right)$ and $\left({\ast}_1\right)$ in $\V\[\U\]$, contradicting \mbox{choice of E}. 
\end{proof}
\begin{Remark}
	If $\A$ is any infinite MAD family in ${\[\omega\]}^{\omega}$ and if ${\E}_0 = \{E \in {\[\omega\]}^{\omega}: \existsinf A \in \A \[\lc E \cap A \rc = \omega \] \}$, then Mathias~\cite{M} showed that ${\E}_0$ is a P-coideal. It is easy to see that for a MAD family in ${\omega}^{\omega}$, ${\E}_0$, as defined in Theorem \ref{thm:IisnotPcoideal}, is not necessarily a P-coideal. This is an interesting difference between the two types of MADness.
\end{Remark}
	Next, we will explore some consequences of Theorem \ref{thm:IisnotPcoideal} for some ideals on $\omega$ that can be naturally defined by using a MAD family of functions $\A \subset {\omega}^{\omega}$.
\begin{Def}
	Let $\A \subset {\omega}^{\omega}$ be a MAD family. We define $\IA = \{a \in \Pset(\omega \times \omega): \exists p \in {\omega}^{a} \forall h \in \A \[\lc p \cap h \rc < \omega \] \}$. Given $E \subset \omega \times \omega$, we define $E(n) = \{m \in \omega: \langle n, m \rangle \in E \}$ and $\dom {(E)} = \{n \in \omega: E(n) \neq 0 \}$.
\end{Def}
\begin{Lemma} \label{lem:Eiscoideal}
	Let $\A \subset {\omega}^{\omega}$ be a MAD family. Let $\E = \{E \in {\[\omega \times \omega\]}^{\omega}: \forall k \in \omega \: \exists a \subset \dom{(E)} \[ \; a \notin \IA\wedge \forall n \in a\lc E(n) \rc > k \] \: \}$. $\I = \Pset(\omega \times \omega) \setminus \E$ is an ideal on $\omega \times \omega$.
\end{Lemma}
\begin{proof}
	It is easy to see that $\I$ is closed under subsets. We will check that it is also closed under unions. Fix $E_0, E_1 \in \I$ and suppose, for a contradiction, that $E_0 \cup E_1 \in \E$. Observe that $\dom{(E_0 \cup E_1)} = \dom{(E_0)} \cup \dom{(E_1)}$ and that for all $n \in \omega$, $\left( E_0 \cup E_1 \right)(n) = E_0(n) \cup E_1(n)$. For each $k \in \omega$ and $i \in \{0, 1 \}$, define $a^i_k = \{n \in \omega: \lc E_i(n) \rc > k \}$. Note that $\dom{(E_i)} = a^i_0 \supset a^i_1 \supset \dotsb$. Since $\IA$ is an ideal, if $a^i_k \in \IA$ for some $k$, then $\forall k' \geq k \[a^i_{k'} \in \IA \]$. Therefore, it follows from our assumption that $E_0$ and $E_1$ are both in $\I$ that for some $k \in \omega$, both $a^0_k$ and $a^1_k$ are in $\IA$. Since $E_0 \cup E_1 \in \E$, $\{n \in \omega: \lc E_0(n) \cup E_1(n) \rc > 2k \} \notin \IA$. Therefore, we may choose $n \notin a^0_k \cup a^1_k$ such that $\lc E_0(n) \cup E_1(n) \rc > 2k$. But since $n \notin a^0_k \cup a^1_k$, $\lc E_0(n) \rc \leq k$ and $\lc E_1 (n)\rc \leq k$, a contradiction.    
\end{proof}

\begin{Theorem} \label{thm:I0Anotpcoideal}
	Let $\A \subset {\omega}^{\omega}$ be a MAD family. If ${\[\omega\]}^{\omega} \setminus \IA$ is a P-coideal, then $\A$ is not analytic. 
\end{Theorem}
\begin{proof}
	Let ${\E}_0$ be defined as in Theorem \ref{thm:IisnotPcoideal} and $\E$ as in Lemma \ref{lem:Eiscoideal}. Let $\I = \Pset(\omega \times \omega) \setminus \E$. Lemma \ref{lem:Eiscoideal} tells us that $\I$ is an ideal. Moreover, if $E \in \E$ and $\{h_0, \dotsc, h_k\} \subset {\omega}^{\omega}$, then there is an infinite partial function $p \subset E$ with $\dom(p) \notin \IA$, which is disjoint from $h_0, \dotsc, h_k$. It follows that there are infinitely many $h \in \A$ such that $\lc E \cap h \rc = \omega$, whence $\E \subset {\E}_0$. Therefore, by Theorem \ref{thm:IisnotPcoideal}, it suffices to show that $\E$ is a P-coideal. Fix a sequence $E_0 \supset E_1 \supset \dotsb$, with $E_i \in \E$. For each $i$ and $k$, define $a^i_k = \{n \in \omega: \lc E_i(n) \rc > k \}$. As before, we have $\dom{(E_i)} = a^i_0 \supset a^i_1 \supset \dotsb$. By assumption, no $a^i_k$ is in $\IA$. We also have $a^0_k \supset a^1_k \supset \dotsb$. Thus, $\langle a^k_k: k \in \omega \rangle$ is a decreasing sequence of sets not in $\IA$. Since we are assuming that ${\[\omega\]}^{\omega} \setminus \IA$ is a P-coideal, there is a set $a \notin \IA$ such that $a \: {\subset}^{\ast} a^k_k$, for all $k$. Let us define a set $E \subset \omega \times \omega$ with $\dom{(E)} = a$ as follows. Let $\langle n_i: i \in \omega \rangle$ enumerate $a$. We may assume that $a \subset a^0_0$. For each $i \in \omega$, let $l_i = \textrm{max} \{k \leq i: n_i \in a^k_k \}$. Note that $n_i \in a^{l_i}_{l_i}$, and hence that $\lc {E}_{l_i}(n_i) \rc > l_i$. Therefore, we may define $E(n_i)$ to be some (arbitrary) subset of ${E}_{l_i}(n_i)$ of size equal to $l_i + 1$. We will check that $E$ is as required.  Since $a \: {\subset}^{\ast} a^k_k$, $\lim l_i = \infty$, and therefore, $\lim \lc E (n_i) \rc = \infty$. As, $\dom{(E)} = a \notin \IA$, this gives us $E \in \E$. Next, we must check that $E \: {\subset}^{\ast} E_k$ for all $k$. Fix $k$. We know that $\forallbutfin i \in \omega \[ l_i  \geq k\]$. Thus $\forallbutfin i \in \omega \[E(n_i) = E_{l_i}(n_i) \subset E_k(n_i) \]$. As each $E(n_i)$ is finite, we get that $E \: {\subset}^{\ast} E_k$. 
\end{proof}
\bibliographystyle{amsplain}
\bibliography{Bibliography}
\end{document}